\documentclass[a4paper,12pt] {amsart}

\usepackage{amssymb}
\usepackage{amsmath}
\usepackage[english,francais]{babel}
\usepackage[active]{srcltx} 

\newtheorem{theorem}{Theorem}[section]

\begin{document}

\title{Formule de Kirillov et conjecture de Guillemin-Sternberg }

\author{Michel Duflo}
\address{Universit\'e Paris Diderot-Paris 7, Institut de Math\'ematiques de Jussieu,
 C.P.~7012\\ 2~place Jussieu,   F-75251 Paris~cedex~05}
\email{duflo@math.jussieu.fr}

\author{ Mich{\`e}le Vergne}
\address{Universit\'e Paris Diderot-Paris 7, Institut de Math\'ematiques de Jussieu,
 C.P.~7012\\ 2~place Jussieu,   F-75251 Paris~cedex~05}
\email{vergne@math.jussieu.fr}

\selectlanguage{english}
\title{Kirillov's formula and Guillemin-Sternberg conjecture}

\begin{abstract}
\selectlanguage{english}

Let $G$ be a connected reductive real Lie group, and $H$ a  compact connected subgroup.
Harish-Chandra associates to a  regular coadjoint  admissible orbit  $M$ of $G$  some unitary representations $\Pi$.
Using the character formula for $\Pi$, we show
that the multiplicities of $\Pi$ restricted to $H$ can be computed  in terms of the fibers of the restriction map  $p:M\to \mathfrak h^*$, if $p$ is proper. In particular, this gives an alternate proof of a result of Paradan.

\vskip 0.5\baselineskip

\selectlanguage{francais}
\noindent{\bf R\'esum\'e} \vskip 0.5\baselineskip \noindent

Soit $G$ un groupe  de Lie r\'eel r\'eductif connexe, et $H$ un sous-groupe   compact connexe.
Harish-Chandra associe  \`a  une orbite coadjointe  r\'eguli\`ere admissible   $M$ de $G$  des repr\'esentations unitaires
$\Pi$.
Gr\^ace aux formules de caract\`ere pour $\Pi$, nous montrons que les multiplicit\'es de la restriction de $\Pi$ \`a $H$ peuvent \^etre calcul\'ees
en fonction des fibres de la projection   $p:M\to \mathfrak h^*$, si $p$ est propre.
En particulier, ceci donne une autre d\'emonstration d'un r\'esultat de Paradan.

\end{abstract}

\maketitle

\selectlanguage{english}

\section{Introduction}

When $M$ is a compact pre-quantizable Hamiltonian manifold for the action of a compact connected Lie group $H$ with moment map $p: M\to {\mathfrak h}^*$, Guillemin and Sternberg  defined  a quantization of $M$,  which is a virtual representation of $H$.  They  proposed formulas for the multiplicities in terms of the reduced ``manifolds" $p^{-1}(v)/H(v)$.
These formulae have been proved in \cite{meinrenken}, and in the close setting of   quantization with $\rho$-correction (also called $Spin_c$-quantization)  in  \cite{paradanspin}. In this note we consider only quantization with  $\rho$-correction.

When $M$ is not compact, it is not clear how to define a quantization of $M$. In the case where $M$ is a
coadjoint admissible orbit, closed and of maximal dimension,
of a real connected reductive Lie group $G$, the representations $\Pi$ associated to $M$ by Harish-Chandra are the natural candidates for the  quantization of $M$.
When $H$ is a maximal compact subgroup of $G$,
Paradan \cite{par1} has shown that  the motto  ``quantization
commutes with reduction" still holds for these non compact Hamiltonian manifolds.
We give another  proof  using   character formulae.
It  holds  for any connected compact subgroup $H$, provided  the moment map $p: M\to {\mathfrak h}^*$ is proper.
However, our proof uses a special feature of these manifolds $M$: their $\hat A$ genus is trivial. So
it  does not extend to representations associated to coadjoint orbits of $G$  which are not of maximal dimension.

\section{Box splines and Dahmen-Micchelli
 deconvolution formula}
Let $V$ be a  finite dimensional real vector space,  $\Lambda\subset V$ a lattice,   and $dv$ the associated Lebesgue measure. For $v\in V$, we denote by $\delta_v$ the $\delta$ measure at $v$, by $\partial_v$ the differentiation in the direction $v$. Let $\Phi:=[\alpha_1,\alpha_2,\ldots, \alpha_N]$ be a list of elements in $\Lambda$ and let $\rho_\Phi:=\frac{1}{2}\sum_{\alpha\in \Phi}\alpha$.
The  centered Box spline $B_c(\Phi)$ is the measure on $V$ such that, for a  continuous function  $F$ on $V$,
$$
\langle B_c(\Phi),F\rangle :=
\int_{-\frac{1}{2}}^{\frac{1}{2}}\cdots \int_{-\frac{1}{2}}^{\frac{1}{2}}
F(\sum_{i=1}^{N} t_i\alpha_i) dt_1\cdots dt_N.
$$
The Fourier transform  ${\hat B}_c(\Phi)(x)$  is the function
$\prod_{\alpha\in \Phi}\frac{e^{i\langle \alpha,x\rangle/2}-e^{-i\langle  \alpha,x\rangle/2}}{i\langle \alpha,x\rangle}$.

Define a series of differential operators on $V$ by
$$
{\hat A}(\Phi):=\prod_{\alpha\in \Phi}
\frac{\partial_\alpha}{e^{\frac{1}{2}\partial_\alpha}-e^{-\frac{1}{2}\partial_\alpha}}
 =1-\frac{1}{24}\sum_{\alpha}(\partial_\alpha)^2+\cdots.
 $$
We assume now that $\Phi$ generates $V$.  A point $v\in V$ is called $\Phi$-regular if $v$ does not lie on any affine hyperplane $\rho_\Phi + \lambda + U$ where $\lambda\in \Lambda$ and $U$ is a hyperplane spanned by elements of $\Phi$.
A connected component ${\mathfrak c}$ of the set of $\Phi$-regular elements is called an alcove.
A vector  $\epsilon \in V$ is called generic if $\epsilon$ does not lie on any   hyperplane $U$.

 Denote by ${\mathcal C}(\Lambda)$ the space of complex valued functions on $\Lambda$.
Let ${\bf h}=(h^{{\mathfrak c}})_{{\mathfrak c}}$ be  a family of polynomial  functions on $V$ indexed by the alcoves, and $\epsilon$ generic.
We define $\lim_{\epsilon}{\bf h}\in \mathcal C(\Lambda)$ by
$\lim_{\epsilon}{\bf h}(\lambda)=h^{{\mathfrak c}}(\lambda)$, where ${\mathfrak c}$ is the alcove  such that $\lambda+t\epsilon\in {\mathfrak c}$ for small  $t>0$.
 If $P$ is a differential operator (or a series of differential operators)  with constant coefficients, we define
 $P({\bf h})=(Ph^{{\mathfrak c}})_{{\mathfrak c}}$.

Consider the Box spline $B_c(\Phi)$. For each alcove ${\mathfrak c}$, there exists a polynomial function $b^{{\mathfrak c}}$ on $V$ such that the measure $B_c(\Phi)$ coincide with $b^{{\mathfrak c}}(v)dv$ on ${\mathfrak c}$. We thus obtain a family ${\bf b}=(b^{{\mathfrak c}})_{{\mathfrak c}}$.
If  $m\in {\mathcal C}(\Lambda)$, the convolution of  the discrete measure $\sum_{\lambda} m(\lambda)\delta_\lambda$ with $B_c(\Phi)$ is given by a polynomial function $b(m)^{{\mathfrak c}}$ on each alcove ${\mathfrak c}$. We denote by ${\bf b}(m)$ the family so obtained.

Recall that the list $\Phi$ is called  unimodular if any basis of $V$ contained in $\Phi$  is a basis of the lattice $\Lambda$.
\begin{theorem} \cite{DM}.
If $\Phi$ is unimodular, then for any $\epsilon$ generic,
$$\lim_{\epsilon}({\hat A}(\Phi){\bf b}(m))=m.$$
\end{theorem}
To explain what happens when  $\Phi$ is not unimodular, we need more notations.
We consider $\Lambda$ as the group of characters of a torus $T$,
and use the notation $s^\lambda$ for the value
of $\lambda \in \Lambda$ at $s\in T$. For $s \in T$,
we denote by $\hat s \in {\mathcal C}(\Lambda)$ the corresponding character of $\Lambda$.
Let ${\mathcal V}(\Phi)$  be the set of $s\in T$ such that the list $\Phi_s:=[\alpha, s^{\alpha}=1]$ still generates $V$ (it is called the vertex set).

Consider a vertex  $s\in {\mathcal V}(\Phi)$ and the convolution product
\begin{equation}\label{eq1}
B_c(s,\Phi):=\left(\prod_{\alpha\in \Phi\setminus \Phi_s} \frac{ \delta_{\alpha/2}-{s}^{-\alpha }\delta_{-\alpha/2}}{1-{s}^{-\alpha }}\right)*B_c(\Phi_s).
\end{equation}

 If $m\in {\mathcal C}(\Lambda)$,  Theorem \ref{deconv} (basically due to Dahmen-Micchelli) below says that we can recover the value  of $m$ at a point   $\lambda_0$  from the knowledge, in a neighborhood of $\lambda_0$, of   the locally polynomial measures (for all $s\in {\mathcal V}(\Phi)$)
\begin{equation}\label{moreimportant}
   b(s,m,\Phi):=\left(\sum_{\lambda} ({\hat s} m)(\lambda)\delta_\lambda\right)*B_c(s,\Phi)
\end{equation}
and moreover a  precise formula is given  in terms of  differential operators.
The measure $ b(s,m,\Phi)$  is given  on each alcove ${\mathfrak c}$ by a polynomial function.
We denote by ${\bf b}(s, m )$ the family of polynomials so obtained.

Define the  series of differential operators
$$E(  s ,\Phi):=
\prod_{\alpha\in \Phi\setminus \Phi_s}\frac{1-s^{-\alpha }}{ e^{\partial_{-\alpha/2}}- s^{-\alpha }e^{{\partial_{\alpha/2}}}}
= 1 + \frac{1}{2}\sum_{\alpha\in \Phi\setminus \Phi_s}  \frac{1+s^{-\alpha}}{1-s^{-\alpha}}\ \partial_{\alpha}+ \cdots
$$
and  $${\hat A}(s,\Phi):=E(s,\Phi){\hat A}(\Phi_s).$$

\begin{theorem}\label{deconv} (see \cite{DCP}).
Let $m\in {\mathcal C}(\Lambda)$.
For any $\epsilon$ generic, we have
$$\sum_{s\in {\mathcal V}(\Phi)} \ {\hat s}^{-1}\lim_\epsilon
{\hat A}(s,\Phi)({\bf b}(s,m))=m.$$
\end{theorem}

\section{Kirillov character formulae}
Let $G$ be a connected reductive real Lie group with Lie algebra $\mathfrak g.$
The function $j_{\mathfrak g}(X):=\det_{\mathfrak g}\left(\frac{e^{ad X/2}-e^{-ad X/2}}{ad X}\right)$ admits a square root $j_{\mathfrak g}^{1/2}(X)$, an analytic function on $\mathfrak g$ with $j_{\mathfrak g}^{1/2}(0)=1$.
 Let $s$ be a semi-simple element of $G$, and ${\mathfrak g}(s)$ its centralizer. The function
 $\det_{{\mathfrak g}/{\mathfrak g(s)}}\left(\frac{1-se^{ad X}}{1-s}\right)$ admits a square root $D^{1/2}(s,X)$, an analytic function on ${\mathfrak g}(s)$ with $D^{1/2}(s,0)=1$.

Let $H$ be a compact connected  group, with Lie algebra $\mathfrak h$.
Let $T$ be a Cartan subgroup of $H$ with Lie algebra $\mathfrak t$. We will apply the results of the previous paragraph to the vector space $V=\mathfrak t^*$ equipped with the lattice $\Lambda\subset {\mathfrak t}^*$  of weights of $T$ (thus $e^{i\lambda}$ is a character of $T$). Let  $W$ be  the Weyl group of $(H,T)$.
Choose a positive system $\Delta^+\subset {\mathfrak t}^*$  for the non zero weights of the adjoint action of $T$ in ${\mathfrak h}_{\mathbb C}$.
For $X\in {\mathfrak t}$, $j_{\mathfrak h}^{1/2}(X)=\prod_{\alpha\in\Delta^+}\frac{e^{i\langle \alpha,X\rangle /2}-e^{-i\langle \alpha,X\rangle /2}}{i\langle \alpha,X\rangle }$.

Let $\rho_H:=\rho_{\Delta^+}$,  and
${\mathfrak t}^*_+$  be the open Weyl  chamber.
   Thus ${\mathfrak t}^*_+$ intersect every orbit of $H$ in ${\mathfrak h}^*$ of maximal dimension in one point.
Consider the set $P_{\mathfrak h}:=(\rho_H+\Lambda) \subset {\mathfrak t}^*$
and $P_{\mathfrak h}^+ :=(\rho_H+\Lambda)\cap {\mathfrak t}^*_+$.
A function  $mult$ on $P_{\mathfrak h}^+$ will be extended to a $W$-anti-invariant function $m$ on $P_{\mathfrak h}$.

The set  $P_{\mathfrak h}^+$
is in one-to-one correspondence $\mu\to \Pi^H(\mu)$  with the dual $\hat H$  of $H$.
The identity $${\rm Tr }\, \Pi^H(\mu)(\exp X)=\sum_{w\in W} \frac{\epsilon(w)e^{\langle iw\mu,X\rangle }}{\prod_{\alpha\in \Delta^+}e^{i\langle \alpha,X\rangle /2}-e^{-i\langle \alpha,X\rangle /2}}$$ holds on ${\mathfrak t}$. This is the Atiyah-Bott fixed-point formula for the index of a twisted Dirac operator on $H\mu$,
so that $\Pi^H(\mu)$ is the quantization $Q(H\mu)$ of the symplectic manifold $H\mu$.

Let $p: M\to {\mathfrak h}^*$ be the  moment map of a connected $H$-hamiltonian manifold $M$. Let $\beta_M$ be the Liouville measure.
The \emph{slice}  $S$  of $M$ is the locally closed  subset $p^{-1}({\mathfrak t}_+^*)$ of $M$.
It is a symplectic submanifold of $M$ with associated Liouville measure $\beta_S$.
If $p$ is proper, the restriction $p^0$ of $p$ to $S$ defines a proper map $S\to {\mathfrak t}^*_+$. We extend the push-forward measure $p^0_*(\beta_S)$
on ${\mathfrak t}_+^*$ to a $W$-anti-invariant signed measure on ${\mathfrak t}^*$ denoted by $DH(M,p)$ (the Duistermaat-Heckman measure). If $S$ is non empty (that is, if $p(M)$ contains an $H$-orbit of maximal dimension), the support of  $DH(M,p)$ is equal to $p(M)\cap \mathfrak t^*$.
Suppose moreover that   there exists regular values  $v\in {\mathfrak t}_+^*$     of $p^0$. At such $v$,
the reduced space $M_v:=p^{-1}(v)/H(v)$ is an orbifold with symplectic form denoted by $\Omega_v$,  and corresponding Liouville measure  $\beta_{M_v}$. By \cite{DH}, the measure  $DH(M,p)$ has a polynomial  density with respect to $dv$ in a neighborhood of $v$, and the value at $v$ is the symplectic volume
$\int_{M_v}e^{\Omega_v/2\pi}
=\int_{M_v} \beta_{M_v}$.

Some  unitary irreducible representations of $G$ can  similarly be associated
to closed admissible orbits of maximal dimension of the coadjoint  representation of $G$. Recall Harish-Chandra parametrization. To simplify, we assume $G$ linear.
Let $f_0 \in {\mathfrak g}^*$ such that ${\mathfrak g}(f_0 )$ (its centralizer in ${\mathfrak g}$) is a Cartan subalgebra  of ${\mathfrak g}$. Denote by ${\tilde G}(f_0 )$  the metaplectic  two fold cover of the stabilizer  $G(f_0 )$ of $f_0 $.
Let  $\tau$ be a   character of
 $\tilde G(f_0 )$ such that $\tau(\exp X)=e^{i\langle f,X\rangle }$ and $\tau(\epsilon)=-1$ if $\epsilon\in {\tilde G}(f_0 )$  projects on $1$ and $\epsilon\neq 1$ (if such a character $\tau$ exists, $f_0$ is called admissible). As explained in \cite{du},   Harish-Chandra associated to this data an irreducible unitary  representation $\Pi^G(f_0 ,\tau)$ of $G$. We consider it as a quantization $Q(M,\tau)$ of $M$. If $f_0$ is admissible and $G(f_0)$ is connected (as is the case when $G(f_0)$ is compact), the character $\tau$ is unique, and  we simply write $Q(M)$ for $Q(M,\tau).$

Irreducible unitary representations of $G$  have a   character, which, by Harish-Chandra theory, is a locally $L^1$ function on $G$. We denote by $\Theta(M,\tau)$ the character of $Q(M,\tau)$.  Similarly, the measure $\beta_M$, considered as a tempered measure on $\mathfrak g^*$,  has a Fourier transform which is a locally $L^1$ function  on $\mathfrak g$.
Kirillov formula (proven in this case by Rossmann \cite{rossmann}) is the equality of locally $L^1$  functions on $\mathfrak g$:
$$
 j_{\mathfrak g}^{1/2}(X)\ \Theta(M,\tau)(\exp X)
 =\int_M e^{i\langle f,X\rangle }d\beta_M(f).
 $$

We suppose that the connected compact group  $H$ is a subgroup of $G$, and we assume that   the projection map $p:M\to \mathfrak h^*$ is proper.
It implies that the restriction $Q(M,\tau)|_H=\sum_{\mu\in P_{\mathfrak h}^+} mult(\mu) \Pi^H(\mu)$ is a sum of irreducible representations of $H$
with finite multiplicities $mult(\mu)$.
 We  associated to $mult(\mu)$  an anti-invariant function  $m(\mu)$ on $P_{\mathfrak h}$, and to the projection $p$ an anti-invariant measure $DH(M,p)$ on ${\mathfrak t}^*$.
Let $\Delta({\mathfrak g}/{\mathfrak h})\subset \mathfrak t^*$ be the list of  weights  for the action of $T$ in ${\mathfrak g}_{\mathbb C}/{\mathfrak h}_{\mathbb C}$.
Choose a  sublist  $\Phi$   so that
$\Delta({\mathfrak g}/{\mathfrak h})$ is the disjoint union of $\Phi$, $-\Phi$ and the zero weights.
 The subsequent definitions do not depend of this choice.
On $\mathfrak t$, we have
$j_{\mathfrak g}^{1/2}(X)=j_{\mathfrak h}^{1/2}(X)
\prod_{\alpha\in\Phi}\frac{e^{i\langle \alpha,X\rangle /2}-e^{-i\langle \alpha,X\rangle /2}}{i\langle \alpha,X\rangle }.$

We assume (and we can easily restrict to this case) that $\mathfrak h$ does not contain any ideal of $\mathfrak g$.
Then the set $\Phi$ generates $\mathfrak t^*$.

Kirillov formula, written  for the characters of $Q(M,\tau)$ and  $Q(H\mu)$, implies the equality of measures on ${\mathfrak t}^*$:
$$(\sum_{\nu\in P_{\mathfrak h}} m(\nu)\delta_\nu)* B_c(\Phi)=DH(M,p).$$
The  measure $DH(M,p)$  is a  polynomial $d^{\mathfrak a}$  on each alcove ${\mathfrak a}={\mathfrak c}+\rho_H$,
 where $\mathfrak c$  is an alcove for the system $\Phi$. For $v\in \mathfrak a$,
   $r(v):=\hat A(\Phi)d^{\mathfrak a}(v)=\int_{M_v}e^{\Omega_v/2\pi} \hat A(M_v)$ where $\hat A(M_v)$ is the $\hat A$ genus of $M_v$.
  This follows from  expressing  the linear variation of $\Omega_v$ in function of the curvature of the principal bundle $(p^0)^{-1}(v)/T$ (\cite{DH}).

In the (rare) case where the system $\Phi$ is unimodular (for example for $G$ the adjoint group of $U(p,q)$, and $H$ the maximal compact subgroup), the orbifold $M_v$ is smooth.
 The value $r(\nu)$ can be defined at any $\nu\in P_{\mathfrak h}$,
 by taking a limit of $r(\nu+t\epsilon)$ where $\nu+t\epsilon$ stay in an alcove $\mathfrak a$, and any $\epsilon$ generic, and coincide with  the number $Q(M_\nu)\in {\mathbb Z}$ defined as the quantization of the (possibly singular)
  reduction $M_\nu$ (\cite{par1}). Thus we obtain
$$Q(M)|_H=\sum_{\nu\in P_{\mathfrak h}^+\cap p(M)} Q(M_\nu)Q(H\nu).$$

We now consider the general case.
 Consider a vertex $s\in T$ for $\Phi$. Let $M(s)\subset M $ be the submanifold  of fixed points of $s$. It may have several connected components.
  It is a symplectic submanifold, and we denote by $\beta_s$ its Liouville measure.
We can define the generalized function $\Theta(f,\tau)(sg)$ where $g\in G$ commutes with $s$. The identity
 \begin{equation}\label{bouaziz}
 j_{\mathfrak g(s)}^{1/2}(X)D^{1/2}(s,X)\Theta(f,\tau)(s\exp X)=
 \int_{M(s)} \epsilon(s,\tau) e^{i\langle f,X\rangle }
  \beta_s\end{equation}
  holds  as an identity of locally $L^1$- functions on $\mathfrak g(s)$ (\cite{Bouaziz}).
Here $\epsilon(s,\tau)$ is a locally constant function on $M(s)$ (\cite{DuVe}).

Recall (\ref{eq1}) the measure $B_c(s,\Phi)$ on $\mathfrak t^*$ associated to $s$.
Denote by $p_s: M(s)\to {\mathfrak h}(s)^*$ the restriction of $p$ to $M(s)$. We define $DH(M,s,\tau)$ as the
sum of the measures $\epsilon(s,\tau)_i DH(M_s^i,p_s^i)$, where $M_s^i$ are the connected components of $M_s$, and $\epsilon(s,\tau)_i$
 the constant value of $\epsilon(s,\tau)$ on $M_s^i$.
The support of $DH(M,s,\tau)$ is contained in the image $p(M)$ of $M$ for any $s$.
Formula (\ref{bouaziz})  implies the identity of measures on ${\mathfrak t}^*$:
$$(\sum_{\nu\in P_{\mathfrak h}} ({\hat s} m)(\nu) \delta_\nu)* B_c(s,\Phi)=DH(M,s,\tau).$$
  Comparing with formula (\ref{moreimportant}), we see that we can compute $m(\nu)$ from the knowledge, in a neighborhood of $\nu$, of Duistermaat-Heckmann measures $DH(M,s,\tau)$  associated to all vertices $s$. In particular $m(\nu)=0$, if $\nu$ is not in the interior of  $p(M)^0$ of $p(M)$. More precisely, Theorem \ref{deconv} and the definition of the quantization of (possibly singular) reduced spaces   gives us
$$Q(M,\tau)|_H=\sum_{\nu\in P_{\mathfrak h}^+\cap p(M)^0} Q(M_\nu,\tau) Q(H\nu).$$

\end{document}